\input amstex
\documentstyle{amsppt}
\magnification=1200
\hsize=6.5true in
\vsize=23.5true cm
\def\version{File Paper137v6a last changed 20.10.00}
\def\nmonth{\ifcase\month\ \or January\or
   February\or March\or April\or May\or June\or July\or August\or
   September\or October\or November\else December\fi}
\NoRunningHeads
\topmatter
\title Algebraic curvature tensors whose skew-symmetric
       curvature operator has constant rank $2$\endtitle
\author Peter Gilkey and Tan Zhang\endauthor
\rightheadtext{Gilkey-Zhang}
\leftheadtext{Heat Trace Asymptotics}
\pageno=1\abstract Let $R$ be an algebraic
curvature tensor for a non-degenerate inner product of signature $(p,q)$ where $q\ge5$. If
$\pi$ is a spacelike $2$ plane, let $R(\pi)$ be the associated skew-symmetric curvature
operator. We classify the algebraic curvature tensors so $R(\cdot)$ has
constant rank $2$ and show these are geometrically realizable by hypersurfaces in flat
spaces. We also classify the Ivanov-Petrova algebraic curvature tensors of rank
$2$; these are the algebraic curvature tensors of constant rank
$2$ such that the complex Jordan normal form of $R(\cdot)$ is constant.\endabstract
\endtopmatter

\newcount\qcts\newcount\qcta\qcta=1\newcount\qcteq\newcount\qcthead
\def\sethead#1{}
\def\setref#1{}
\def\seteqn#1{}
\def\eqnpbg{\ifnum\qcteq=1 a\else\ifnum\qcteq=2 b\else
   \ifnum\qcteq=3 c\else\ifnum\qcteq=4 d\else
   \ifnum\qcteq=5 e\else\ifnum\qcteq=6 f\else\ifnum\qcteq=7 g\else
   \ifnum\qcteq=8 h\else\ifnum\qcteq=9 i\else
   \ifnum\qcteq=10 j\else\ifnum\qcteq=11 k\else\ifnum\qcteq=12 l\else
   \ifnum\qcteq=13 m\else\ifnum\qcteq=14 n\else\ifnum\qcteq=15 o\else
   \ifnum\qcteq=16 p\else\ifnum\qcteq=17 q\else
   \ifnum\qcteq=18 r\else\ifnum\qcteq=19 s\else\ifnum\qcteq=20 t\else
   \ifnum\qcteq=21u\else\ifnum\qcteq=22 v\else
   \ifnum\qcteq=23 w\else\ifnum\qcteq=24 x\else
   \ifnum\qcteq=25 y\else\ifnum\qcteq=26 z\else *
    \fi\fi\fi\fi\fi\fi\fi\fi\fi\fi\fi\fi\fi\fi\fi\fi\fi\fi\fi\fi\fi\fi\fi\fi\fi\fi}
\newcount\qct\newcount\qcta\qct=0\qcta=1
\def\pbgkey#1{}
\font\pbglie=eufm10 
\def\so{\text{\pbglie so}}
\def\t{\Cal{T}}\def\a{\text{\pbglie{a}}}\def\A{\Cal{A}}
\def\Pb{\text{\pbglie{b}}}
\def\range{\operatorname{Range}}
\def\rank{\operatorname{Rank}}
\def\span{\operatorname{Span}}
\def\spec(#1){\operatorname{Spec(#1)}}
\def\Gr{\operatorname{Gr}}
 \def\qctSA{1}
 \def\arefa{1.1}
 \def\arefaa{1.1.a}
 \def\arefab{1.1.b}
 \def\arefac{1.1.c}
 \def\arefad{1.1.d}
 \def\arefae{1.1.e}
 \def\arefb{1.2}
 \def\arefc{1.3}
 \def\arefca{1.3.a}
 \def\arefd{1.4}
 \def\arefdA{1.4.a}
 \def\arefe{1.5}
 \def\arefda{1.5.a}
 \def\qctSB{2}
 \def\brefa{2.1}
 \def\brefb{2.2}
 \def\brefc{2.3}
 \def\brefd{2.4}
 \def\brefe{2.5}
 \def\breff{2.6}
 \def\breffa{2.6.a}
 \def\brefg{2.7}
 \def\brefh{2.8}
 \def\qctSC{3}
 \def\crefa{3.1}
 \def\crefaa{3.1.a}
 \def\crefab{3.1.b}
 \def\crefac{3.1.c}
 \def\crefb{3.2}
 \def\crefba{3.2.a}
 \def\crefbb{3.2.b}
 \def\qctSD{4}
 \def\drefa{4.1}
 \def\drefb{4.2}
 \def\drefc{4.3}
 \def\qctSE{5}
 \def\erefa{5.1}
 \def\erefaa{5.1.a}
 \def\erefb{5.2}
 \def\qctSF{6}
 \def\qctSG{7}
 \def\grefa{7.1}
 \def\qctSH{8}
 \def\hrefa{8.1}
 \def\hrefb{8.2}
 \def\hrefc{8.3}
\def\refBBG{1}
\def\refBBGZ{2}
\def\refChia{3}
\def\refGra{4}
\def\refGVV{5}
\def\refGia{6}
\def\refGib{7}
\def\refGilw{8}
\def\refGis{9}
\def\refGLS{10}
\def\refGSem{11}
\def\refIP{12}
\def\refKo{13}
\def\refOss{14}
\def\refZ{15}
\sethead\qctSA\head\S\qctSA\ Introduction\endhead
\setref\arefa
Let ${}^g\nabla$ be the Levi-Civita connection of a pseudo-Riemannian
manifold $(M,g)$ of signature $(p,q)$. Let
${}^gR(x,y):={}^g\nabla_x{}^g\nabla_y-{}^g\nabla_y{}^g\nabla_x-{}^g\nabla_{[x,y]}$ be the
Riemann curvature operator and let ${}^gR(x,y,z,w)$ be the associated curvature tensor. We
have:
\seteqn\arefaa
\seteqn\arefab
\seteqn\arefac
\seteqn\arefad
$$\leqalignno{
&g({}^gR(x,y)z,w)={}^gR(x,y,z,w),&(\text{\arefaa})\cr
&{}^gR(x,y,z,w)=-{}^gR(y,x,z,w)=-{}^gR(x,y,w,z),&(\text{\arefab})\cr
&{}^gR(x,y,z,w)={}^gR(z,w,x,y),\text{ and}&(\text{\arefac})\cr
&{}^gR(x,y,z,w)+{}^gR(y,z,x,w)+{}^gR(z,x,y,w)=0.&(\text{\arefad})\cr}$$

It is convenient to work in a purely algebraic context. Let $\langle\cdot,\cdot\rangle$ be
a non-degenerate symmetric inner product of signature $(p,q)$ on a finite dimensional real
vector space $V$. We say that a $4$ tensor
$R\in\otimes^4(V^*)$ is an {\it algebraic curvature tensor} if $R$ satisfies the identities
of equations (\arefab), (\arefac), and (\arefad); the {\it associated algebraic curvature
operator} $R(\cdot,\cdot)$ is then defined by equation (\arefaa). Given an inner product 
$\langle\cdot,\cdot\rangle_P$ and an algebraic
curvature tensor $R_P$ on the tangent space $T_PM$ at a point $P\in M$, there is the germ of
a metric $g$ on $M$ so $g|_{T_PM}=\langle\cdot,\cdot\rangle_P$ and so ${}^gR|_{T_PM}=R_P$.
Thus algebraic curvature tensors are important in differential geometry.

A central problem in differential geometry is to understand the relationship between
algebraic properties of the curvature tensor ${}^gR$ and the underlying geometry of the
manifold. The full curvature tensor is in general quite difficult to work with so one often
uses the curvature tensor to define a natural endomorphism of the tangent bundle. One
wants to know the geometric consequences that follow if such an operator is assumed to have
constant eigenvalues or to have constant rank. The Jacobi operator, the Stanilov operator,
the Szab\'o operator, and the skew-symmetric curvature operator are examples of such natural
operators. We refer to
\cite{\refGia, \refGis} for further details; in the interests of brevity
we shall content ourselves here with a brief discussion of the Jacobi and skew-symmetric
curvature operators as motivation. 

We say that a vector $v$ is {\it spacelike} if
$\langle v,v\rangle>0$. We say that a subspace $\pi$ of $V$ is {\it spacelike} if the
restriction of
$\langle\cdot,\cdot\rangle$ to $\pi$ is positive definite. Let
$$\eqalign{
   &V^+:=\{v\in V:\langle v,v\rangle>0\},\ 
 \Bbb{S}^+(V):=\{v:\langle v,v\rangle=1\},\text{ and }
  \Gr_{(0,2)}^+(V,\langle\cdot,\cdot\rangle)}$$
be the set of spacelike vectors, the set of unit spacelike
vectors, and the set of oriented spacelike $2$ planes in $V$. Let $R$ be an algebraic
curvature tensor. We define the {\it Jacobi operator}
$J_R(x):y\rightarrow R(y,x)x$. We say that
$R$ is {\it Osserman} if the complex eigenvalues of $J_R(\cdot)$ are constant on
$\Bbb{S}^+(V)$. If
$(M,g)$ is a local rank $1$ symmetric space or is flat, then the local isometries act
transitively on the bundle $\Bbb{S}^+(TM)$ and hence the eigenvalues of $J_R$ are constant on
$\Bbb{S}^+(TM)$. In the Riemannian setting ($p=0$), Osserman
\cite{\refOss} wondered if the converse held; this question has become known as the Osserman
conjecture. Chi \cite{\refChia} established the Osserman conjecture for Riemannian manifolds
of dimension $m$ for $m\equiv1$ mod $2$, for $m\equiv2$ mod $4$, and for $m=4$. The
corresponding conjecture for algebraic curvature tensors fails if
$m\equiv0$ mod $4$ \cite{\refGib}. In the Lorentzian setting ($p=1$),
the Osserman conjecture has been established \cite{\refBBG, \refGra}. The Osserman
conjecture is false in the higher signature setting ($p>1$) \cite{\refBBGZ,
\refGVV}.

Let $R$ be an algebraic curvature tensor. If $\{v_1,v_2\}$ is an oriented basis for a
spacelike
$2$ plane
$\pi$, the {\it skew-symmetric curvature operator} is defined independently of the basis by:
\seteqn\arefae
$$R(\pi):=\{\langle v_1,v_1\rangle\langle v_2,v_2\rangle-\langle
v_1,v_2\rangle^2\}^{-1/2} R(v_1,v_2).\tag\arefae$$
We say $R$ has {\it rank $r$} if
$\rank(R(\pi))=r$ for any oriented spacelike $2$ plane $\pi$. The following
theorem was proved using topological methods \cite{\refGLS, \refZ}. It shows that $r=2$ in
many cases.
\proclaim{\arefa\ Theorem}
Let
$R$ be an algebraic curvature tensor of rank $r$.\roster
\smallskip\item Let $p\le1$. Let $q=5$, $q=6$, or
$q\ge9$. Then $r=2$.
\smallskip\item Let $p=2$. Let $q\ge10$. Assume neither $q$ nor $q+2$ are
powers of $2$. Then
$r=2$.
\endroster
\endproclaim

If $(p,q)=(0,4)$, if $(p,q)=(2,2)$, or if $(p,q)=(4,0)$, then the
algebraic curvature tensors of rank $4$ form an open non-empty subset of the space
of all algebraic curvature tensors so the $4$ dimensional setting is exceptional. It is not
known if there exist algebraic curvature tensors of rank greater than $2$ for other
values of
$(p,q)$. Kowalski et al \cite{\refKo} has shown in the Riemannian setting that an algebraic
curvature tensor of rank $4$ in dimension $4$ must have both positive and negative sectional
curvatures.

In light of Theorem \arefa, we shall focus our attention on the algebraic curvature tensors
of rank $2$. If $\phi$ is a linear map from $V$ to $V$, we define a $4$ tensor
$R_\phi$ and associated operator by twisting the tensor of constant sectional curvature $+1$:
$$\eqalign{&R_\phi(x,y)z:=\langle\phi y,z\rangle\phi x
   -\langle\phi x,z\rangle\phi y,\text{ and}\cr
    &R_\phi(x,y,z,w):=\langle\phi y,z\rangle\langle\phi x,w\rangle
   -\langle\phi x,z\rangle\langle\phi y,w\rangle.}$$
The adjoint is defined by the identity
$\langle\phi x,y\rangle=\langle x,\phi^*y\rangle$. We say that
$\phi$ is {\it admissible} if $\phi^*=\phi$ and if the kernel of $\phi$ contains no
spacelike vectors. We can classify rank $2$ algebraic curvature tensors:
\setref\arefb
\proclaim{\arefb\ Theorem} Let $q\ge5$.
\roster
\smallskip\item If $\phi$ is admissible, then
$\pm R_\phi$ is a rank $2$ algebraic curvature tensor.
\smallskip\item If $R$ is a rank $2$ algebraic curvature tensor, then $R=\pm
R_\phi$ for some admissible $\phi$.
\endroster\endproclaim

We say that an algebraic curvature tensor $R$ of rank $r$ is {\it Ivanov-Petrova} if the
complex Jordan normal form of $R(\pi)$ is the same for every oriented spacelike $2$ plane
$\pi$; these tensors were first classified for $m=4$ and
$p=0$ by Ivanov and Petrova \cite{\refIP}. The following result generalizes results of
\cite{\refGLS, \refZ} from the Riemannian and Lorentzian settings to
arbitrary signatures.

\setref\arefc
\proclaim{\arefc\ Theorem} Let $q\ge5$. A rank $2$ algebraic curvature tensor $R$ is
Ivanov-Petrova if and only if $R=\pm R_\phi$ for some admissible $\phi$ where
one of the following conditions holds:
\roster
\smallskip\item There exists $C\ne0$ so $\langle\phi v_1,\phi v_2\rangle=C\langle
v_1,v_2\rangle$ for all $v_i\in V$.
\smallskip\item The range of $\phi$ is totally isotropic, i.e.
$\langle\phi v_1,\phi v_2\rangle=0$ for all $v_i\in V$.
\endroster\endproclaim

To prove Theorem \arefb, it is convenient to decouple the domain and range and to establish a
slightly more general classification result. Let $\langle\cdot,\cdot\rangle_A$ and
$\langle\cdot,\cdot\rangle_B$ be non-degenerate symmetric inner products of signatures
$(p_A,q_A)$ and $(p_B,q_B)$ on $A$ and
$B$. Let
$$\so(B,\langle\cdot,\cdot\rangle_B):=\{\phi:\phi+\phi^*=0\}$$
be the set of all skew-symmetric linear maps; this is the Lie algebra of the special
orthogonal group defined by
$\langle\cdot,\cdot\rangle_B$.  Let $T$ be an alternating bilinear map from
$A\otimes A$ to $\so(B,\langle\cdot,\cdot\rangle_B)$. The associated $4$ tensor
$T(a_1,a_2,b_1,b_2)$ defined by equation (\arefaa) has the symmetries of equation
(\arefab); conversely given a $4$ tensor with the symmetries of equation (\arefab), we can
use equation (\arefaa) to define an alternating map from $A\otimes A$ to
$\so(B,\langle\cdot,\cdot\rangle_B)$. We use equation (\arefae) to extend
$T$ to a map from
$\Gr_{(0,2)}^+(A,\langle\cdot,\cdot\rangle_A)$ to
$\so(B,\langle\cdot,\cdot\rangle_B)$ and suppose this extension has constant rank $2$.
Such a map will be said to be {\it admissible}.

Let $\phi$ be a linear map from $A$ to $B$, let $\xi\in B$, and let $\chi$ be an
alternating bilinear map from $A\otimes A$ to $B$. Define:
\seteqn\arefca
$$\eqalign{
     &T_{\chi,\xi}(a_1,a_2)b:=\langle\chi(a_1,a_2),b\rangle_B\xi
     -\langle\xi,b\rangle_B\chi(a_1,a_2),\cr
&T_{\chi,\xi}(a_1,a_2,b_1,b_2):=\langle\chi(a_1,a_2),b_1\rangle_B\langle\xi,b_2\rangle_B
    -\langle\xi,b_1\rangle_B\langle\chi(a_1,a_2),b_2\rangle_B,\cr
     &T_\phi(a_1,a_2)b:=\langle\phi(a_2),b\rangle_B\phi(a_1)
     -\langle\phi(a_1),b\rangle_B\phi(a_2),\text{ and}\cr
&T_\phi(a_1,a_2,b_1,b_2)=\langle\phi(a_2),b_1\rangle_B\langle\phi(a_1),b_2\rangle_B
                    -\langle\phi(a_1),b_1\rangle_B\langle\phi(a_2),b_2\rangle_B.}
    \tag\arefca$$
Since the symmetries of equation (\arefab) are satisfied,
$T_{\chi,\xi}$ and $T_\phi$ are alternating bilinear maps from $A\otimes A$ to
$\so(B,\langle\cdot,\cdot\rangle_B)$.
\setref\arefd
\proclaim{\arefd\ Theorem} Let $q_A\ge5$.\roster
\smallskip\item We have the following criteria for admissibility:
\smallskip\itemitem{\rm 1-a)} $T_\phi$ is admissible if and
only if $\ker(\phi)\cap A^+=\emptyset$.
\smallskip\itemitem{\rm 1-b)} $T_{\xi,\chi}$ is admissible if and only if
$\{\chi(a_1,a_2),\xi\}$ are linearly independent vectors whenever $\{a_1,a_2\}$ spans a
spacelike $2$ plane.
\smallskip\item Let $T$ be admissible. Then either $T=T_{\chi,\xi}$ or $T=\pm T_\phi$ for
suitable $(\chi,\xi)$ or $\phi$.
\endroster
\endproclaim

Let $R$ be an algebraic curvature tensor. We say that a pseudo-Riemannian manifold
$(M,g)$ is a {\it geometric realization} of $R$ at a point $P\in M$ if there is an
isomorphism $\psi$ from $T_PM$ to $V$ so that $\psi^*\{\langle\cdot,\cdot\rangle_V\}=g(P)$
and so that $\psi^*\{R\}={}^gR(P)$. It is a classical result that every algebraic curvature
tensor has a geometric realization.  Let $R$ be a rank $2$ algebraic curvature tensor. By
Theorem
\arefb, there exists an admissible $\phi$ so $R=\varepsilon R_\phi$ where
$\varepsilon=\pm1$. Let $e_{m+1}\cdot\Bbb{R}$ be a $1$ dimensional vector space. Define:
\seteqn\arefdA
$$\eqalign{
  &W:=V\oplus e_{m+1}\cdot\Bbb{R},\ 
   F_\phi(v):=v\oplus\textstyle\frac12\langle\phi(v),v\rangle e_{m+1},\text{ and}\cr
  &\langle v_1\oplus r_1e_{m+1},v_2\oplus r_2e_{m+1}\rangle_W=\langle
v_1,v_2\rangle_V+\varepsilon
   r_1r_2.}\tag\arefdA$$ 
Then $F_\phi$ defines the germ of an embedding of
$V$ into $W$ as a hypersurface which is defined in some small neighborhood of the origin. Let
$g_\phi$ be the associated first fundamental form. We take the canonical identification of
$T_0V=V$. We say that $(M,g)$ has rank $r$ if ${}^gR$ has rank
$r$ at each point of the manifold. 

\setref\arefe
\proclaim{\arefe\ Theorem} Let $R=\pm R_\phi$ be a rank $2$ algebraic curvature tensor on
$(V,\langle\cdot,\cdot\rangle_V)$. Then $(V,g_\phi)$ is the germ of a rank
$2$ pseudo-Riemannian manifold which realizes $R$ at $P=0$.
\endproclaim

We say that a pseudo-Riemannian manifold $(M,g)$ is {\it Ivanov-Petrova} if
$(M,g)$ has constant rank $r$ and if ${}^gR$ is Ivanov-Petrova at each point of
$M$; the complex Jordan normal form can vary with the point in question but the
rank is assumed constant. The Riemannian Ivanov-Petrova manifolds have been
classified if $m\ge4$ and $m\ne7$ \cite{\refGilw, \refGLS, \refIP}; the case
$m=7$ is exceptional and the classification is not complete although some partial results
are known \cite{\refGSem}. We refer to \cite{\refZ} for some partial results in
the pseudo-Riemannian setting. This classification shows that there exist
Ivanov-Petrova algebraic curvature tensors which are not geometrically realized
by Ivanov-Petrova pseudo-Riemannian manifolds.

Here is a brief outline to the paper. In \S\qctSB, we establish some technical results. In
Lemma \brefa, we discuss linearizations of a projective map. In Lemmas \brefb, \brefc\ and
\brefd, we establish some of the elementary properties of
$\so(B,\langle\cdot,\cdot\rangle_B)$. In Lemma \breff, we give necessary and sufficient
conditions for $T_\phi$ to be an algebraic curvature tensor if
$\ker(\phi)$ contains no spacelike vector. We conclude \S\qctSB\ by proving Theorems
\arefb\ (1) and \arefd\ (1).

Let $\Bbb{P}(A)$ and $\Bbb{P}(B)$ be the projective spaces defined by $A$ and $B$. Let
$[a]=a\cdot\Bbb{R}\in\Bbb{P}(A)$ and $[b]:=b\cdot\Bbb{R}\in\Bbb{P}(B)$ be the lines spanned
by $0\ne a\in A$ and $0\ne b\in B$. Let $T$ be admissible. We define
\seteqn\arefda
$$\Phi([a]):=\cap_{\pi\in \Gr_{(0,2)}^+(A):[a]\subset\pi}\range(T(\pi)).\tag\arefda$$
In \S\qctSC, we assume $p_A=p_B=0$. We show in Lemma \crefa\ that $\dim\Phi=1$. If $\Phi$ is
constant, we choose $0\ne\xi\in\Phi$, we set $\chi(a_1,a_2):=T(a_1,a_2)\xi$, and we show
$T=CT_{\chi,\xi}$. If $\Phi$ is non-constant, we show there
exists a map $\phi$ which linearizes $\Phi$ and we show $T=CT_\phi$ for some constant
$C$. This completes the proof of Theorem \arefd\ (2) in the
Riemannian setting.

In \S\qctSD, we assume $p_B=0$ and let $p_A$ be arbitrary. Let $S$ be a maximal
spacelike subspace of $A$ and let $T_S=T|_{S\otimes S}$. We use the results of
\S\qctSC\ to express $T_S=T_{\chi(S),\xi(S)}$ or $T_S=T_{\phi(S)}$. The main technical difficulty is then to extend these tensors to
all of
$A$ and thereby prove Theorem
\arefd\ (2) in this setting. In
\S\qctSE, we use the map
$\psi$ defined in Lemma \brefb\ (1) to show that the signature of the metric on $B$ plays no
role and to complete the proof of Theorem \arefd\ (2). We then use Theorem \arefd\ (2) to
complete the proof of Theorem \arefb\ (2). In \S\qctSF, we use Theorem \arefb\ to prove
Theorem
\arefc. 

In \S\qctSG, we review some facts concerning the geometry of hypersurfaces and use
Theorem \arefb\ to prove Theorem
\arefe. In
\S\qctSH, we construct some examples. In Example
\hrefa, we show that the assumption
$q_A\ge5$ in Theorem \arefd\ is essential by constructing an admissible map $T$ from
$\Bbb{R}^{(p,4)}\otimes\Bbb{R}^{(p,4)}$ to $\so(4)$ so that $T\ne\pm T_\phi$ and
$T\ne T_{\chi,\xi}$ for any $\phi$ or
$(\chi,\xi)$. In Example \hrefb, we give an algebraic
curvature tensor which has constant eigenvalues and constant rank $2$ on the set of
oriented spacelike $2$ planes but which is not Ivanov-Petrova; this shows that it is
necessary to consider the complex Jordan normal form and not simply the eigenspace structure
and rank. A subspace
$\pi$ is said to be {\it timelike} if the restriction of
$\langle\cdot,\cdot\rangle$ to $\pi$ is negative definite. By replacing
$\langle\cdot,\cdot\rangle$ by $-\langle\cdot,\cdot\rangle$ we can interchange the roles of
timelike and spacelike subspaces. We use equation (\arefae) to define
$R(\pi)$ for oriented timelike $2$ planes. In Example \hrefc, we exhibit a rank $2$
Ivanov-Petrova algebraic curvature tensor which does not have
constant rank on the set of timelike $2$ planes. Thus the behavior for spacelike and
timelike planes can be different.

It is a pleasant task to acknowledge helpful conversations with Dr. Inis Kath dealing with
stylistic matters in the paper.

\sethead\qctSB\head\S\qctSB\ Some technical preliminaries\endhead
Let $\Bbb{P}(A^+)\subseteq\Bbb{P}(A)$ be the subset of spacelike lines in $A$. Let
$\Phi$ be a map from $\Bbb{P}(A^+)$ to $\Bbb{P}(B)$. We say that a linear map $\phi$
from $A$ to $B$ {\it linearizes }$\Phi$ if $\ker(\phi)$ contains no
spacelike vectors and if
$[\phi(a)]=\Phi[a]$ for all $[a]\in\Bbb{P}(A^+)$. Note that vectors $\{a_1,...,a_p\}$ are
linearly independent if and only if $a_1\wedge...\wedge a_p\ne0$ in $\Lambda^p(A)$.
\setref\brefa
\proclaim{\brefa\ Lemma} Let $\Phi:\Bbb{P}(A^+)\rightarrow\Bbb{P}(B)$.\roster
\smallskip\item If $q_A\ge3$ and if $\phi_i$ are linearizations of $\Phi$, then $\phi_1$ is a
multiple of $\phi_2$.
\smallskip\item Let $p_A=0$ and let $q_A\ge5$. Let $\Cal{F}$ be a family of codimension $1$
subspaces of $A$. Suppose $\Phi$ is linearizable on every $F\in\Cal{F}$ and that every $2$
plane is contained in some element of $\Cal{F}$. Then $\Phi$ is linearizable.
\endroster\endproclaim

\demo{Proof} The spacelike vectors form an open subset of $A$ which generates $A$
additively. Thus to prove
$\phi_1=c\phi_2$ it suffices to prove that there exists $c$ so $\phi_1(a)=c\phi_2(a)$ for
every spacelike $a$. Let $a_1$ and $a_2$ be spacelike. Since
$q_A\ge3$, we can choose $a_0$ spacelike so $a_0\perp a_1$ and $a_0\perp a_2$.  Since
$[\phi_1a_i]=[\phi_2 a_i]$, there are non-zero real numbers $c_i$ so that
$\phi_1(a_i)=c_i\phi_2(a_i)$. To prove assertion (1) it suffices to show $c_1=c_0=c_2$.
Since the roles of $a_1$ and
$a_2$ are symmetric, we must only show $c_0=c_1$. Since
$\phi_2$ is injective on the spacelike $2$ plane $\span\{a_0,a_1\}$,
$0\ne\phi_2(a_0)\wedge\phi_2(a_1)$. Choose
$c_{01}$ so that
$\phi_1(a_0+a_1)=c_{01}\phi_2(a_0+a_1)$. We show $c_0=c_{01}=c_1$ and complete the
proof of assertion (1) by
computing:
$$\eqalign{
   &c_{01}\phi_2(a_0)+c_{01}\phi_2(a_1)=c_{01}\phi_2(a_0+a_1)=\phi_1(a_0+a_1)\cr
  &\qquad=\phi_1(a_0)+\phi_1(a_1)
  =c_0\phi_2(a_0)+c_1\phi_2(a_1).\cr}$$

Let $\phi_F$ be the linearization of $\Phi$ on $F\in\Cal{F}$. Fix a basepoint $\bar
F\in\Cal{F}$. Both $\phi_{F}$ and $\phi_{\bar F}$ are linearizations of
$\Phi$ on $F\cap\bar F\ne\{0\}$. Since $\dim(F\cap\bar F)\ge q_A-2\ge3$, we can use assertion
(1) to see $\phi_{F}=c(F)\phi_{\bar F}$ on $F\cap\bar F$. We replace $\phi_F$ by
$c(F)^{-1}\phi_F$ to assume without loss of generality that $\phi_F=\phi_{\bar F}$ on
$F\cap\bar F$ for all $F\in\Cal{F}$. Let $F_i\in\Cal{F}$. Since
$\phi_{F_1}=\phi_{\bar F}=\phi_{F_2}$ on
$F_1\cap F_2\cap\bar F\ne\{0\}$, we use assertion (1) to see $\phi_{F_1}=\phi_{F_2}$ on
the (possibly) larger intersection $F_1\cap F_2$. Thus $\phi(a)=\phi_F(a)$ is well
defined on $\cup_{F\in\Cal{F}}F=A$. Since every $2$ plane is contained in some element
of $\Cal{F}$, $\phi$ is linear and linearizes $\Phi$.
\qed\enddemo

We can relate $\langle\cdot,\cdot\rangle_B$ to a positive definite inner product.

\setref\brefb
\proclaim{\brefb\ Lemma} There exists a self-adjoint linear map $\psi$ from $B$ to $B$ with
$\psi^2=id$ so that
$(b_1,b_2)_+:=\langle\psi b_1,b_2\rangle_B$ is a positive definite symmetric inner
product on $B$, so that
$\langle b,\tilde b\rangle_B=(\psi b,\tilde b)_+=(b,\psi\tilde b)_+$, and so that
$(b,\tilde b)_+=\langle\psi b,\tilde b\rangle_B=\langle b,\psi\tilde b\rangle_B$.
The map $\t\rightarrow\psi\t$ is an isomorphism between
$\so(B,\langle\cdot,\cdot\rangle_B)$ and $\so(B,(\cdot,\cdot)_+)$.
\endproclaim

\demo{Proof} Choose a basis $\{\Pb_i\}$ for $B$ so
$\langle\Pb_i,\Pb_j\rangle_B=0$ for $i\ne j$, so
$\langle\Pb_i,\Pb_i\rangle_B=-1$ for $i\le p_B$, and so 
$\langle\Pb_i,\Pb_i\rangle_B=1$ for $i>p_B$. Let $\{\Pb^i\}$ be the corresponding dual basis
for $B^*$. Define:
$$(b_1,b_2)_+:=\textstyle\sum_i\Pb^i(b_1)\Pb^i(b_2)\text{ and }
  \psi(b_1):=-\textstyle\sum_{i\le p_B}\Pb^i(b_1)\Pb_i
        +\textstyle\sum_{i>p_B}\Pb^i(b_1)\Pb_i.$$
The first identities now follow. We complete the proof by establishing the following chain
of equivalent statements:
\medbreak\qquad\qquad (1) $\t\in\so(B,\langle\cdot,\cdot\rangle_B)$.
\smallbreak\qquad\qquad (2) $\langle\t b,\tilde b\rangle_B+\langle b,\t\tilde b\rangle_B=0$
for all
$b,\tilde b\in B$.
\smallbreak\qquad\qquad (3) $(\psi\t b,\tilde b)_++(b,\psi\t\tilde b)_+=0$ for all $b,\tilde
b\in B$.
\smallbreak\qquad\qquad (4) $\psi\t\in\so(B,(\cdot,\cdot)_+)$. \qed\enddemo

If $T\in\so(B,\langle\cdot,\cdot\rangle_B)$, define $\omega(T)\in\Lambda^2(B^*)$ by
$\omega(T)(b_1,b_2):=\langle Tb_1,b_2\rangle_B$. The correspondence $T\rightarrow\omega(T)$
identifies $\so(B,\langle\cdot,\cdot\rangle_B)$ with $\Lambda^2(B^*)$.

\setref\brefc
\proclaim{\brefc\ Lemma} Let $0\ne\t\in\so(B,\langle\cdot,\cdot\rangle_B)$. Then $\rank(\t)$
is even. We have
$\rank(\t)=2$ if and only if $\omega(\t)\wedge\omega(\t)=0$.\endproclaim

\demo{Proof}
Since $\rank(\t)=\rank(\psi\t)$ and 
$\omega_{\langle\cdot,\cdot\rangle_B}(\t)=\omega_{(\cdot,\cdot)_+}(\psi\t)$, we may suppose
the metric is positive definite. Choose
an orthonormal basis $\{\Pb_i\}$ for $B$ and $\lambda_\mu>0$ so
$$\eqalign{
   &\t\Pb_{j}=0\text{ for }j>2\ell\text{ and }
   \t\Pb_{2\mu}=\lambda_\mu\Pb_{2\mu-1},\ 
    \t\Pb_{2\mu-1}=-\lambda_\mu\Pb_{2\mu}\text{ for }\mu\le\ell.\cr
   }$$
Thus $\rank(\t)=2\ell$ is even. Let $\{\Pb_i^*\}$ be the dual basis for $B^*$. We use the
following decomposition to complete the proof:
$$\omega(\t)=\textstyle\sum_{1\le\mu\le\ell}
   \lambda_\mu\Pb_{2\mu-1}^*\wedge\Pb_{2\mu}^*.\ \qed$$\enddemo

\setref\brefd
\proclaim{\brefd\ Lemma} Let $\t_i\in\so(B,\langle\cdot,\cdot\rangle_B)$. Assume
$\rank(\t_i)=2$.\roster

\smallskip\item Then
$\t_1$ is a multiple of $\t_2$ if and only if
$\range(\t_1)=\range(\t_2)$.
\smallskip\item Let $p_B=0$. Let $0\ne b\in\range(\t_1)$.
Then $\{b,\t_1b\}$ is an orthogonal basis for
$\range(\t_1)$.
\endroster\endproclaim

\demo{Proof} We use Lemma \brefb\ to assume the metric is positive definite in the proof of
assertion (1). Since
$\rank(\t_i)=2$, $\t_i$ is a multiple of a 90 degree rotation in the $2$ plane
$\range(\t_i)$ and vanishes on $\range(\t_i)^\perp$. Assertions (1) and (2) now follow.
\qed
\enddemo

The following technical observation will prove useful in later sections.

\setref\brefe
\proclaim{\brefe\ Lemma} Let $T$ be a bilinear map from $A\otimes A$ to a vector space $C$.
Assume that $T(a_1,a_2)=0$ whenever $a_1$ and $a_2$ span a spacelike $2$ plane. Then
$T=0$.\endproclaim

\demo{Proof} The set of such tensors $a_1\otimes a_2$ generates $A\otimes A$
additively. Since $T$ is bilinear, we may conclude $T=0$. \qed\enddemo

We continue our preparatory steps with a final technical Lemma.

\setref\breff
\proclaim{\breff\ Lemma} Let $\phi$ be a linear map from $A$ to $B$ so
$\ker(\phi)$ contains no spacelike vector. Let $q_A\ge3$.\roster
\smallskip\item If $T_{\tilde\phi}=\pm T_\phi$, then $\tilde\phi=\pm\phi$.
\smallskip\item Let $A=B=V$. Then
$T_\phi$ is an algebraic curvature tensor if and only if $\phi=\phi^*$.
\endroster\endproclaim

\demo{Proof} Use equation (\arefda) to define $\Phi$. Let $a_1$ be spacelike.
Since
$q_A\ge3$, we can choose $a_2$ and $a_3$ so $\span\{a_1,a_2,a_3\}$ is spacelike. Then
$0\ne\phi(a_1)\wedge\phi(a_2)\wedge\phi(a_3)$ so
$$\eqalign{
\phi(a_1)&\cdot\Bbb{R}
\subseteq\Phi([a_1])\subseteq\range(T_\phi(a_1,a_2))\cap\range(T_\phi(a_1,a_3))\cr
&\subseteq\span\{\phi(a_1),\phi(a_3)\}\cap\span\{\phi(a_1),\phi(a_3)\}
   \subseteq\phi(a_1)\cdot\Bbb{R}.}$$
Thus $\Phi([a_1])=\phi(a_1)\cdot\Bbb{R}$ and $\phi$ linearizes $\Phi$. If $T_\phi=\pm
T_{\tilde\phi}$, then $\tilde\phi$ also linearizes $\Phi=\tilde\Phi$ so by Lemma \brefa\ (1),
$\tilde\phi=c\phi$. As $T_{c\phi}=c^2T_\phi$, $c^2=\pm1$ so $c=\pm1$ and
assertion (1) follows.
 
Let $A=B=V$. Since $T_\phi$ defines an alternating bilinear map from $V\otimes V$ to
$\so(V,\langle\cdot,\cdot\rangle)$, the curvature symmetry (\arefab) is immediate. We
compute:
$$\eqalign{
   &T_\phi(a_3,a_4,a_1,a_2)=
     \langle\phi(a_4),a_1\rangle\langle\phi(a_3),a_2\rangle-
     \langle\phi(a_3),a_1\rangle\langle\phi(a_4),a_2\rangle\cr
  &\qquad=\langle\phi^*(a_1),a_4\rangle\langle\phi^*(a_2),a_3\rangle-
     \langle\phi^*(a_1),a_3\rangle\langle\phi^*(a_2),a_4\rangle\cr
  &\qquad=T_{\phi^*}(a_1,a_2,a_3,a_4).}$$
Thus symmetry of equation (\arefac) is satisfied if and only if $T_\phi=T_{\phi^*}$ or
equivalently using assertion (1) if and only if $\phi=\pm\phi^*$. 
If $\phi^*=\phi$, we verify that the Bianchi identities of equation (\arefad) are satisfied
by computing:
$$\eqalignno{
 &T_\phi(x,y)z+T_\phi(y,z)x+T_\phi(z,x)y\cr
=&\langle\phi y,z\rangle\phi x-\langle\phi x,z\rangle\phi y
  +\langle\phi z,x\rangle\phi y-\langle\phi y,x\rangle\phi z
  +\langle\phi x,y\rangle\phi z-\langle\phi z,y\rangle\phi x\cr
=&(\langle\phi y,z\rangle-\langle\phi z,y\rangle)\phi x
  +(\langle\phi z,x\rangle-\langle\phi x,z\rangle)\phi y
  +(\langle\phi x,y\rangle-\langle\phi y,x\rangle)\phi z
  =0.\cr}$$
Conversely if the Bianchi identities are satisfied, we suppose that $\phi=-\phi^*$ and
argue for a contradiction. We compute:
\seteqn\breffa
$$\eqalign{
0=&T_\phi(a_1,a_2)a_3+T_\phi(a_2,a_3)a_1+T_\phi(a_3,a_1)a_2\cr
    =&\langle\phi a_2,a_3\rangle\phi a_1+\langle\phi a_3,a_1\rangle\phi a_2
    +\langle\phi a_1,a_3\rangle\phi a_3\cr
    &-\langle\phi a_1,a_3\rangle\phi a_2-\langle\phi a_2,a_1\rangle\phi a_3
    -\langle\phi a_3,a_2\rangle\phi a_1\cr
     =&2\langle\phi a_2,a_3\rangle\phi a_1+2\langle\phi a_3,a_1\rangle\phi a_2
    +2\langle\phi a_1,a_2\rangle\phi a_3.\cr}\tag\breffa$$
Fix $a_1$ spacelike. Then $\phi a_1\ne 0$. Since $q_A\ge3$, we may choose $a_2$ so
$a_2\perp\phi a_1$ and so
$\{a_1,a_2\}$ is a spacelike orthonormal set. If $a_3\perp\phi a_1$, since $\phi a_1\ne0$,
we use equation (\breffa) to see that $a_3\perp\phi a_2$. Thus $\phi a_2$ is a multiple of
$\phi a_1$. This is not possible as $\phi$ is injective on spacelike subspaces. Consequently
we have
$\phi=\phi^*$. \qed\enddemo

\setref\brefg
\subhead\brefg\ The proof of Theorem \arefd\ (1)\endsubhead We have
$\rank(T_\phi)\le2$ and
$\rank(T_{\chi,\xi})\le2$. By Lemma \brefc, the rank is even. Thus
$T_\phi$ has constant rank $2$ if and only if $\{\phi(a_1),\phi(a_2)\}$ are linearly
independent vectors whenever $\{a_1,a_2\}$ spans a spacelike $2$ plane. This is equivalent
to the condition that $\ker(\phi)$ contains no spacelike vector. Similarly
$T_{\chi,\xi}$ has constant rank $2$ if and only if $0\ne\chi(a_1,a_2)\wedge\xi$ whenever
$\{a_1,a_2\}$ spans a spacelike $2$ plane. \qed

\setref\brefh
\subhead\brefh\ The Proof of Theorem \arefb\ (1)\endsubhead Suppose that $\phi=\phi^*$ and
that
$\ker(\phi)$ contains no spacelike vectors. By Theorem \arefd\ (1), $R_\phi$ has constant
rank $2$. By Lemma \breff\ (2), $R_\phi$ is an algebraic curvature tensor. \qed

\sethead\qctSC\head\S\qctSC\ Properties of the map $\Phi$\endhead 
We suppose $p_A=p_B=0$ and work in the Riemannian setting in this section. The following
result concerning $\Phi$ is central to our investigation.
\setref\crefa
\proclaim{\crefa\ Lemma} Let $p_A=p_B=0$ and let
$T$ be admissible. Let $0\ne a_1\wedge a_2\wedge a_3$.
\roster
\smallskip\item Let $q_A\ge3$. We have
$\dim\{\range(T(a_1,a_2))\cap\range(T(a_1,a_3))\}=1$.
\smallskip\item Let $q_A\ge5$. Let $\Phi$ be defined by equation {\rm (\arefda)}. We have
$\dim(\Phi([a_1]))=1$.
\smallskip\item Let $q_A\ge5$. We have $\Phi([a_1])=\range(T(a_1,a_2))
\cap\range(T(a_1,a_3))$.
\smallskip\item Let $q_A\ge5$. The map $\Phi:\Bbb{P}(A)\rightarrow\Bbb{P}(B)$ is either
injective or constant.
\endroster\endproclaim

\demo{Proof} We adapt arguments used to study IP algebraic curvature
tensors in \cite{\refGLS}. Let $\pi_{ij}:=\span\{a_i,a_j\}$ and
let $u:=\dim\{\range(T(\pi_{12}))+\range(T(\pi_{13}))\}$. Then
$$\dim\{\range(T(\pi_{12}))\cap\range(T(\pi_{13}))\}=4-u.$$
We may estimate that $4\ge u\ge\dim\range(T(\pi_{12}))=2$. To prove assertion (1) we must
show $u=3$ i.e. that $u\ne2$ and $u\ne 4$. If $u=2$, then
$\range(T(\pi_{12}))=\range(T(\pi_{13}))$ is a $2$ dimensional space so
$T(\pi_{12})=cT(\pi_{13})$ by Lemma \brefd\ (1). Thus
$T(a_1,a_2-ca_3)=0$ which is false. This shows that $u\ne2$. Let $\{\Pb_i\}$ be an
orthonormal basis for
$B$ so
\seteqn\crefaa
$$\eqalign{
   &\range(T(\pi_{12}))=\span\{\Pb_1,\Pb_2\},\ 
    T(\pi_{13})\Pb_1\in\span\{\Pb_1,\Pb_2,\Pb_3\},\text{ and}\cr
  &\range(T(\pi_{12}))+\range(T(\pi_{13}))
    \subseteq\span\{\Pb_1,\Pb_2,\Pb_3,\Pb_4\}.}\tag\crefaa$$
Define $\t_{ij}\in\so(B,\langle\cdot,\cdot\rangle_B)$ with
$\omega(\t_{ij})=\Pb_i^*\wedge\Pb_j^*$ by
$\t_{ij}:\Pb_i\rightarrow\Pb_j$ ,$\Pb_j\rightarrow-\Pb_i$, and $\Pb_k\rightarrow
0$ for $k\ne i,j$.
By rescaling $T$, we may assume
$T(\pi_{12})=\t_{12}$. We use equation (\crefaa) to expand:
$$\eqalign{
  &T(\pi_{12})=\t_{12},\ \omega(T(\pi_{12}))=\Pb_1^*\wedge\Pb_2^*,\cr
  &T(\pi_{13})=c_{12}\t_{12}+c_{13}\t_{13}+c_{23}\t_{23}+c_{24}\t_{24}+c_{34}
\t_{34},\text{ and}\cr
   &\omega(T(\pi_{13}))=c_{12}\Pb_1^*\wedge\Pb_2^*
   +c_{13}\Pb_1^*\wedge\Pb_3^*
    +c_{23}\Pb_2^*\wedge\Pb_3^*+c_{24}\Pb_2^*\wedge\Pb_4^*+c_{34}\Pb_3^*\wedge\Pb_4^*.\cr}
$$
Let $\t(\varepsilon):=T(a_1,a_2+\varepsilon a_3)$.
Since $\rank(\t(\varepsilon))=2$, we use Lemma \brefc\ to
see:
\seteqn\crefab
$$0=\omega(\t(\varepsilon))\wedge\omega(\t(\varepsilon))=
     \{2\varepsilon c_{34}+\varepsilon^2(c_{12}c_{34}-c_{13}c_{24})\}
     \Pb_1^*\wedge \Pb_2^*\wedge\Pb_3^*\wedge \Pb_4^*.\tag\crefab$$ 
Since equation (\crefab) holds for all $\varepsilon\in\Bbb{R}$, $c_{34}=0$ and
$c_{13}c_{24}=0$. We complete the proof of assertion (1) by dealing with the two possible
cases:
\roster
\medbreak\item If
$c_{24}=0$, then
$T(\pi_{13})=c_{12}\t_{12}+c_{13}\t_{13}+c_{23}\t_{23}$ so\newline
$\range(T(\pi_{13}))\subseteq\span\{e_1,e_2,e_3\}$ and
$u\le3$. 
\smallskip\item If
$c_{13}=0$, then
$T(\pi_{13})=c_{12}\t_{12}+c_{23}\t_{23}+c_{24}\t_{24}$ so\newline
$\range(T(\pi_{13}))\subseteq\span\{e_1,e_2,c_{23}e_3+c_{24}e_4\}$
and $u\le3$.\endroster

\medbreak\noindent Suppose assertion (2) fails; we argue for a contradiction. Choose $a_4$
with $0\ne a_1\wedge a_4$ so
\seteqn\crefac
$$\cap_{2\le i\le 4}\range(T(a_1,a_i))=\{0\}.\tag\crefac$$
If $2\le i<j\le4$, then
$\range(T(a_1,a_i))\ne\range(T(a_1,a_j))$ so
$0\ne a_1\wedge a_i\wedge a_j$. We use assertion (1) to define:
$$\eqalign{
   &L_{ij}:=\range(T(a_1,a_i))\cap\range(T(a_1,a_j))\in\Bbb{P}(B)\cr
   &E:=L_{23}+L_{24}+L_{34}.}$$
Let $\{i,j,k\}$ be a permutation of $\{2,3,4\}$. Then
$$\{0\}=\range(T(a_1,a_i))\cap\range(T(a_1,a_j))\cap\range(T(a_1,a_k))
   =L_{ij}\cap L_{ik}.$$
Thus $L_{ij}$ and $L_{ik}$ are distinct lines which are contained in $\range(T(a_1,a_i))$.
This shows:
$$\range(T(a_1,a_i))=L_{ij}\oplus L_{ik}\text{ so }\range(T(a_1,a_i))\subseteq E.$$
Let $[a_5]\perp[a_1]$. We wish to show $\range(T(a_1,a_5))\subseteq E$. If 
$0=a_1\wedge a_i\wedge a_5$ for some $i$ with $2\le i\le 4$, then
$T(a_1,a_5)=cT(a_1,a_i)$ so $\range(T(a_1,a_5))\subseteq E$.
Thus we suppose that 
$0\ne a_1\wedge a_i\wedge a_5$ for all
$2\le i\le 4$ and use assertion (1) to define the lines
$$L_{i5}:=\range(T(a_1,a_i))\cap\range(T(a_1,a_5))\in\Bbb{P}(B)\text{ for }2\le i\le 4.$$
If
$L_{25}=L_{35}=L_{45}$, then
$L_{25}=L_{i5}\subseteq\range(T(a_1,a_i))$ for $2\le i\le 4$
which contradicts equation (\crefac). Thus at least two of the lines $L_{i5}$ are
distinct lines which are contained in $\range(T(a_1,a_5))$. Consequently the lines
$\{L_{25},L_{35},L_{45}\}$ span
$\range(T(a_1,a_5))$ and 
$$\range(T(a_1,a_5))\subseteq E.$$
This shows that the map $\psi:a_4\rightarrow\range(T(a_1,a_4))$ is a well defined map from
the projective space $\Bbb{P}(a_1^\perp)$ to the unoriented Grassmannian
$\operatorname{Gr}_2(E)$ of
$2$ planes in $E$. Suppose given distinct lines $[a_4]$ and $[a_5]$ in $\Bbb{P}(a_1^\perp)$.
Since $0\ne a_1\wedge a_4\wedge a_5$, we use assertion (1) to see
$\psi([a_4])\ne\psi([a_5])$. Thus $\psi$ is injective. Since $\dim(E)\le3$, $\dim
\operatorname{Gr}_2(E)\le 2$. It is immediate from the definition that $\psi$ is continuous.
We can therefore use invariance of domain to see $\dim(\Bbb{P}(a_1^\perp))\le2$. Since
$\dim(\Bbb{P}(a_1^\perp))=q_A-2\ge3$, this contradiction completes
the proof of assertion (2). Assertions (3) now follow. In \S\qctSH, we shall present
an example that shows assertion (2) fails if $q_A=4$.

Suppose that $\Phi:\Bbb{P}(A)\rightarrow\Bbb{P}(B)$ is not injective. Choose distinct lines
$[a_i]$ so that $\Phi([a_1])=\Phi([a_2])$. If
$[a_3]\in\Bbb{P}(A)-\Bbb{P}(\span\{a_1,a_2\})$, then
$$\Phi([a_1])=\Phi([a_2])\subseteq\range(T(a_3,a_1))\cap\range(T(a_3,a_2))=\Phi([a_3]).$$
Thus $\Phi$ is constant on $\Bbb{P}(A)-\Bbb{P}(\span\{a_1,a_2\})$. We use assertion (2) to
see $\Phi$ is continuous. It now follows that $\Phi$ is constant on all of $\Bbb{P}(A)$.
\qed\enddemo 

Theorem \arefd\ (2) follows for positive definite metrics from the following Lemma:
\setref\crefb
\proclaim{\crefb\ Lemma} Let $p_A=p_B=0$ and let $q_A\ge5$. Let $T$ be admissible.\roster
\smallskip\item Let $\Phi$ be constant. Fix $\xi\in\Bbb{S}^+(\Phi)$.
Set $\chi(a_1,a_2):=-T(a_1,a_2)\xi$. Then $T=T_{\chi,\xi}$.
\smallskip\item Let $\Phi$ be non-constant. Then there exists a linearization $\phi$ of
$\Phi$ so $T=\pm T_\phi$.
\endroster\endproclaim

\demo{Proof} Suppose $\Phi:\Bbb{P}(A)\rightarrow\Bbb{P}(B)$ is constant. Choose
$\xi\in\Bbb{S}^+(\Phi)$. Then
$\xi\in\Phi[a]$ for all $a$ and we can define an alternating bilinear map $\chi$ from
$A\otimes A$ to $\xi^\perp\subseteq B$ by setting
$\chi(a_1,a_2):=-T(a_1,a_2)\xi\in B$. If $0\ne a_1\wedge a_2$, then we have
$$\range(T(a_1,a_2))=\span\{\xi,\chi(a_1,a_2)\}=\range(T_{\chi,\xi}(a_1,a_2)).$$
Thus by Lemma \brefd\ (1), $T(a_1,a_2)=cT_{\chi,\xi}(a_1,a_2)$. Since
$$T(a_1,a_2)\xi=-\chi(a_1,a_2)=T_{\chi,\xi}(a_1,a_2)\xi,$$
 $c=1$ so
$T(a_1,a_2)=T_{\chi,\xi}(a_1,a_2)$ if $0\ne a_1\wedge a_2$. Lemma \brefe\ now shows
$T=T_{\chi,\xi}$ which completes the proof of assertion (1).

\medbreak  If $\Phi$ is not constant, then $\Phi$ is injective by Lemma \crefa\ (4). The
unit spheres
$\Bbb{S}^+(A)$ and $\Bbb{S}^+(B)$ are the universal
covers of $\Bbb{P}(A)$ and $\Bbb{P}(B)$. We lift $\Phi$ to define:
$$\Phi^+:\Bbb{S}^+(A)\rightarrow \Bbb{S}^+(B)\text{ so }[\Phi^+(a)]=\Phi([a]).$$
Let $\{\a_i\}$ be a basis for $A$
and let $\{\a^i\}$ be the corresponding dual basis for $A^*$. Let $i<j$. Since $\Phi([\a_i])$
and
$\Phi([\a_j])$ are distinct lines which are contained in 
$T(\a_i,\a_j)$, 
$$\range(T(\a_i,\a_j))=\span\{\Phi^+(\a_i),\Phi^+(\a_j)\}.$$
Since $T$ is bilinear, we have $T(a_1,a_2)=\sum_{i,j}\a^i(a_1)\a^j(a_2)T(\a_i,\a_j)$ so:
$$\eqalign{
   &\range(T(a_1,a_2))\subseteq\textstyle\sum_{i,j}\range(T(\a_i,\a_j))
     \subseteq\span\{\Phi^+(\a_1),...,\Phi^+(\a_{q_A})\}.\cr}$$ 
Thus replacing $B$ by $\span\{\Phi^+(\a_1),...,\Phi^+(\a_{q_A})\}$ we
may assume without loss of generality
$$\dim(B)\le\dim(A).$$

Since $\Bbb{P}(A)$ is compact, $\range(\Phi)$ is compact and hence closed. We have that
$\Phi$ is injective and that $\dim(\Bbb{P}(B))\le\dim(\Bbb{P}(A))$. Thus we use invariance of
domain to see that
$\dim(\Bbb{P}(B))=\dim(\Bbb{P}(A))$ and that $\Phi$ is an open map. Thus
$\range(\Phi)$ is an open subset of $\Bbb{P}(B)$. Since $\Bbb{P}(B)$ is connected and
$\range(\Phi)$ is non-empty, $\Phi$ is surjective. Since $\Bbb{P}(A)$ is compact and since
$\Bbb{P}(B)$ is Hausdorff, the fundamental theorem of point set topology shows that $\Phi$
is a homeomorphism from $\Bbb{P}(A)$ onto $\Bbb{P}(B)$. It now follows that the lift
$\Phi^+$ is a homeomorphism from $\Bbb{S}^+(A)$ onto $\Bbb{S}^+(B)$. Since $\Phi$ induces an
isomorphism from the fundamental group of $\Bbb{P}(A)$ to the fundamental group of
$\Bbb{P}(B)$, we have:
$$\Phi^+(-a)=-\Phi^+(a).$$

We now show that $\Phi$ is linearizable; such a linearization is unique up to
constant multiple by Lemma \brefa\ (1). Let $0\ne a_1\wedge a_2$. We first show
$\Phi([a_1+a_2])\subseteq\Phi([a_1])+\Phi([a_2])$. Since $q_B=q_A\ge5$ and since $\Phi^+$ is
surjective, we choose unit vectors $a_3$ and $a_4$ so that
\seteqn\crefba
$$\Phi^+(a_3)\perp\Phi^+(a_4)\text{ and
}\Phi^+(a_3)\perp\{\Phi^+(a_1),\Phi^+(a_2),\Phi^+(a_1+a_2)\}\text{ for }i=3,4.
\tag\crefba$$
Let $i=3,4$. Then $\range(T(a_i,a_1+a_2))\subseteq
\range(T(a_i,a_1))+\range(T(a_i,a_2))$ so we have
$\Phi^+(a_1+a_2)\in\span\{\Phi^+(a_i),\Phi^+(a_1),\Phi^+(a_2)\}$.
We now use equation (\crefba) to see
\seteqn\crefbb
$$\Phi^+(a_1+a_2)\in\span\{\Phi^+(a_1),\Phi^+(a_2)\}.\tag\crefbb$$

We define
$$\A([a]):=\{a_1\in A-\{0\}:\Phi([a_1])\perp\Phi([a])\}\cup\{0\}.$$
It is immediate that $\A([a])$ is closed under scalar
multiplication. Let $a_i\in\A([a])$. We wish to show that
$\a_1+\a_2\in\A([a])$. This is immediate if $0=a_1\wedge a_2$; if $0\ne a_1\wedge a_2$, then
we use equation (\crefbb) to see $a_1+a_2\in\A([a])$. Thus $\A([a])$ is a linear subspace of
$A$. Since $\Phi$ is a homeomorphism from $\Bbb{P}(\A([a]))$ to $\Bbb{P}(\Phi^+(a)^\perp)$,
we use invariance of domain to see
$$\dim\Bbb{P}(\A([a]))=\dim\Bbb{P}(\Phi^+(a)^\perp)\text{ so }
  \dim(\A([a]))=q_A-1.$$
We apply Lemma \brefa\ (2) to the family of codimension $1$ subspaces $\A([a])$. 
Set 
$$\phi_a(a_1):=T(a,a_1)\Phi^+(a)\text{ for }a_1\in\A([a]).$$
Since $T$ is
bilinear, $\phi_a$ is linear. If $a_1\ne0$, then
$\Phi([a_1])\perp\Phi([a])$ so $0\ne a\wedge a_1$. By Lemma \brefd\ (2),
$\{\Phi^+([a]),\phi_a(a_1)\}$ is an orthogonal basis for $\range(T(a,a_1))$.
Since $\Phi([a])\perp\Phi([a_1])$, $\phi_a(a_1)\in\Phi([a_1])$ so $\phi_a$ is a
linearization of $\Phi$ on $\A([a])$. Let $a_i\in A$. Since $q_A\ge3$ and since $\Phi$ is
bijective, we may choose $a\in A$ so $\Phi([a])\perp\Phi([a_1])$ and
$\Phi([a])\perp\Phi([a_2])$. Thus $a_i\in\A([a])$ so every $2$ plane is contained in some
element of $\Cal{F}$. Thus by Lemma \brefa\ (2) there exists a linearization $\phi$ of
$\Phi$ on all of $A$.

We use equation (\arefca) to define $T_\phi$. Let $0\ne a_1\wedge
a_2$. Because
$$T_\phi(a_1,a_2)\in\so(B,\langle\cdot,\cdot\rangle_B)\text{ and }
    \range(T_\phi)=\range(T)=\span\{\phi(a_1),\phi(a_2)\},$$
we may apply Lemma \brefa\ (4) to see
$T_\phi(a_1,a_2)=\mu(a_1,a_2)T(a_1,a_2)$. If $0\ne a_1\wedge a_2\wedge a_3$, then $T(a_1,a_2)$ and
$T(a_1,a_3)$ are linearly independent maps by
Lemma \crefa\ (1). We compute:
$$\eqalign{
T_\phi(a_1,a_2+a_3)=&\mu(a_1,a_2+a_3)T(a_1,a_2+a_3)\cr
=&\mu(a_1,a_2+a_3)T(a_1,a_2)+\mu(a_1,a_2+a_3)T(a_1,a_3)\cr
=&T_\phi(a_1,a_2)+T_\phi(a_1,a_3)\cr
=&\mu(a_1,a_2)T(a_1,a_2)+\mu(a_1,a_3)T(a_1,a_3).\cr}$$
This implies that $\mu(a_1,a_2)=\mu(a_1,a_2+a_3)=\mu(a_1,a_3)$ as
$\{T(a_1,a_2),T(a_1,a_3)\}$ is a linearly independent subset of
$\so(B,\langle\cdot,\cdot\rangle_B)$. Consequently
$\mu(a_1,a_2)$ is independent of $a_2$; a similar argument shows it is independent of $a_1$.
We denote this common value by $\mu$ and express
$T(a_1,a_2)=\mu T_\phi(a_1,a_2)$ if $0\ne a_1\wedge a_2$.
We use Lemma \brefe\ to see $T=\mu T_\phi$ and replace
$\phi$ by $\sqrt{|\mu|}\phi$ to complete the proof. \qed\enddemo

\sethead\qctSD\head\S\qctSD\ Extending linearizations\endhead

Let $\Cal{S}_r$ be the set of all spacelike subspaces of $A$ which have dimension at least
$r$. We have the following inclusions:
$\Cal{S}_{q_A}\subsetneq...\subsetneq\Cal{S}_3$. Let $S\in\Cal{S}_3$ and $\tilde
S\in\Cal{S}_3$.  We say that $\{S_i\}_{i\in\Bbb{Z}}$ is a
{\it chain linking $S$ and $\tilde S$} if $S=S_i$ for some $i$, if $\tilde S=S_j$ for
some $j$, if $S\cap\tilde S\subset S_i\cap S_{i+1}$ for all $i$, and if
$S_i\cap S_{i+1}\in\Cal{S}_3$ for all $i$. We omit the proof of the
following Lemma in the interests of brevity as the proof is straightforward.

\setref\drefa
\proclaim{\drefa\ Lemma} Let $q_A\ge5$. Any two elements of $\Cal{S}_3$ can be linked by a
chain.
\endproclaim

If $S\in\Cal{S}_3$, let
$T_S$ be the restriction of $T$ to $S\otimes S$ and let
$$\Phi_S([a]):=\cap_{[\tilde a]\in\Bbb{P}(S),\ [\tilde a]\ne[a]}\range(T_S(a,\tilde a)).$$

\setref\drefb
\proclaim{\drefb\ Lemma} Let $p_B=0$, let $q_A\ge5$, and let $T$ be admissible.\roster
\smallskip\item If $S_1\in\Cal{S}_3$, then $\dim\Phi_{S_1}([a])=1$ for all
$[a]\in\Bbb{P}(S_1)$.
\smallskip\item If $S_1\in\Cal{S}_3$, if $S_2\in\Cal{S}_3$, and if $S_1\subseteq S_2$, then
$\Phi_{S_2}|_{\Bbb{P}(S_1)}=\Phi_{S_1}$.
\smallskip\item There exists $\Phi:\Bbb{P}(A^+)\rightarrow\Bbb{P}(B)$ so
$\Phi|_{\Bbb{P}(S)}=\Phi_S$ for all $S\in\Cal{S}_3$.
\smallskip\item Either the map $\Phi$ is constant or all the maps $\Phi_S$ are injective.
\endroster\endproclaim

\demo{Proof} Let $S_1\in\Cal{S}_3$, let $S_2\in\Cal{S}_3$, and let $S_1\subseteq S_2$.
Choose $S_3\in\Cal{S}_{q_A}$ so $S_2\subseteq S_3$. Let $[a]\in\Bbb{P}(S_1)$. We have
$\Phi_{S_3}([a])\subseteq\Phi_{S_2}([a])\subseteq\Phi_{S_1}([a])$ so by
Lemma \crefa\ (1,2):
$$1=\dim\{\Phi_{S_3}([a])\}\le\dim\{\Phi_{S_2}([a])\}\le\dim\{\Phi_{S_1}([a])\}=1.$$
Assertions (1) and (2) follow. Fix $[a]\in\Bbb{P}(A^+)$. Let $S$ and $\tilde S$ be any
two spaces in $\Cal{S}_3$ which contain $[a]$. We use Lemma \drefa\ to link $S$ and $\tilde
S$ by a chain $\{S_i\}$. By assertion (2), 
$$\Phi_{S_i}([a])=\Phi_{S_i\cap
S_{i+1}}([a])=\Phi_{S_{i+1}}([a])\text{ for all }i.$$
Consequently $\Phi_S([a])=\Phi_{\tilde S}([a])$ and $\Phi([a]):=\phi_S([a])$ is well defined
and independent of the particular $S$ chosen. This proves assertion (3).

We use Lemma \crefa\ (4) to prove assertion (4). Suppose there exists $S\in\Cal{S}_3$ so
$\Phi_S$ is not injective. Choose $S_1\in\Cal{S}_{q_A}$ so $S\subset S_1$. Then $\Phi_{S_1}$
is not injective so $\Phi_{S_1}$ and hence $\Phi_S$ is constant by Lemma \crefb\ (4).
Use Lemma \drefb\ to link $S$ to any other element $\tilde S\in\Cal{S}_3$ by a chain
$\{S_i\}$. If
$\Phi_{S_i}$ is not injective, then
$\Phi_{S_i}$ is constant. Thus $\Phi_{S_i\cap S_j}$ is constant and hence $\Phi_{S_{i+1}}$
is constant and takes the same values. Thus $\Phi_S$ is not injective implies $\Phi_{\tilde
S}$ and $\Phi_S$ are both constant and take the same value. \qed\enddemo

We can now generalize Lemma \crefb\ to the case $p_A\ne0$:

\setref\drefc
\proclaim{\drefc\ Lemma} Let $p_B=0$ and let $q_A\ge5$. Let $T$ be admissible.\roster
\smallskip\item Let $\Phi$ be constant. Fix $\xi\in\Bbb{S}^+(\Phi)$.
Set $\chi(a_1,a_2):=-T(a_1,a_2)\xi$. Then $T=T_{\chi,\xi}$.
\smallskip\item Let $\Phi$ be non-constant. Then there exists a linearization $\phi$ of
$\Phi$ so $T=\pm T_\phi$.
\endroster\endproclaim

\demo{Proof} Suppose that $\Phi$ is constant.
If $a_1$ and $a_2$ span a spacelike $2$ plane, then Lemma \crefb\ (2) shows
that $T(a_1,a_2)=T_{\chi,\xi}(a_1,a_2)$. Lemma \brefe\ now implies $T=T_{\chi,\xi}$ which
proves assertion (1).

Suppose that $\Phi$ is non-constant. We use Lemma \drefb\ (4) to see $\Phi_S$ is
injective for all $S\in\Cal{S}_3$.
As $\Bbb{S}^+(A)$ is the universal cover of $\Bbb{P}(A^+)$, we lift
$\Phi$ to define
$$\Phi^+:\Bbb{S}^+(A)\rightarrow \Bbb{S}^+(B)
  \text{ so }[\Phi^+(a)]=\Phi([a])\ \forall\ a\in
  \Bbb{S}^+(A).$$
Let $S\in\Cal{S}_3$. Choose $\bar S\in\Cal{S}_{q_A}$ so $S\subset\bar S$. Since $\Phi_{\bar
S}$ is injective, we apply Lemma \crefb\ (2) to $T_S$ to find $\phi_{\bar S}$ so
$T_{\bar S}=T_{\phi_{\bar S}}$; we let $\phi_S=\phi_{\bar S}|_S$. Then $T_S=T_{\phi_S}$. By
Lemma \breff\ (2), $\phi_S$ is determined up to sign. We use $\Phi^+$ to normalize the choice
of sign by requiring
$$(\Phi^+(a),\phi_S(a))_B>0\text{ for all }a\in\Bbb{S}^+(S).$$
If $S_1\in\Cal{S}_3$, if $S_2\in\Cal{S}_3$, and if $S_1\subset S_2$, we have
$\phi_{S_2}|_{S_1}=\phi_{S_1}$. Thus we may use Lemma \drefa\ to find
$$\phi:A^+\rightarrow B\text{ so }\phi|_S=\phi_S\text{ for all }S\in\Cal{S}_3.$$

We now extend $\phi$ to all of $A$. Fix $a\in A$. Let $a_1\in\Bbb{S}^+(A)$ and
$a_2\in\Bbb{S}^+(A)$ satisfy $a\perp a_i$. Choose
$a_3\in\Bbb{S}^+(A)$ so
$\{a,a_1,a_3\}$ and
$\{a,a_2,a_3\}$ are orthogonal sets. Choose
$\varepsilon_i\in\Bbb{R}$ so that
$|\varepsilon_i^2|>|\langle a,a_i\rangle_A|$ for $1\le i\le 3$. As $a+\varepsilon_ia_i\in
A^+$, $\phi(a+\varepsilon_ia_i)-\phi(\varepsilon_ia_i)$ is defined. As
$$\eqalign{
  &\span\{a+\varepsilon_1a_1,a_2,\varepsilon_3a_3\}\in\Cal{S}_3,
  \ \span\{\varepsilon_1\alpha_1,\alpha_2,\varepsilon_3a_3\}\in\Cal{S}_3,
   \cr
  &\span\{a+\varepsilon_3a_3,\varepsilon_1a_1,a_2\}\in\Cal{S}_3,}$$
and as $\phi_S$ is linear if $S\in\Cal{S}_3$, we may compute:
$$\eqalign{
   &\phi(a+\varepsilon_1a_1)-\phi(\varepsilon_1a_1)\cr
  =&\phi(a+\varepsilon_1a_1+\varepsilon_3a_3)-\phi(\varepsilon_3a_3)
    -\phi(\varepsilon_1a_1+\varepsilon_3a_3)+\phi(\varepsilon_3a_3)\cr
  =&\phi(a+\varepsilon_1a_1+\varepsilon_3a_3)-\phi(\varepsilon_1a_1)
     -\phi(\varepsilon_1a_1+\varepsilon_3a_3)+\phi(\varepsilon_1a_1)\cr
  =&\phi(a+\varepsilon_3a_3)-\phi(\varepsilon_3a_3).}$$
As the roles of $a_1$ and $a_2$ are symmetric, this shows that
$$\phi(a+\varepsilon_1a_1)-\phi(\varepsilon_1a_1)
   =\phi(a+\varepsilon_3a_3)-\phi(\varepsilon_3a_3)
   =\phi(a+\varepsilon_2a_2)-\phi(\varepsilon_2a_2)$$
is independent of the choice of $(a_1,\varepsilon_1)$. If $a\in A^+$, then this
difference yields $\phi(a)$. Consequently we may extend $\phi$ from $A^+$ to $A$ by
defining
$$\phi(a):=\phi(a+\varepsilon_1a_1)-\phi(\varepsilon_1a_1).$$
It is immediate that $\phi(\varrho a)=\varrho\phi(a)$. To check $\phi$ is linear, we must
show $\phi$ is additive. If
$\{a,\tilde a\}$ are given, choose $a_1\in\Bbb{S}^+(A)$ and $\tilde a_1\in\Bbb{S}^+(A)$ so
$\{a,a_1,\tilde a_1\}$ and $\{\tilde a,a_1,\tilde a_1\}$ are orthonormal sets. Since
$\{a+\tilde a,a_1+\tilde a_1\}$ is an orthogonal set and $a_1+\tilde a_1\in A^+$, we may
compute
$$\eqalign{
   \phi(a+\tilde a)=&\phi(a+\tilde a+\varepsilon(a_1+\tilde a_1))
     -\varepsilon\phi(a_1+\tilde a_1)\cr
   =&\phi(a+\varepsilon a_1)+\phi(\tilde a+\tilde\varepsilon a_1)
    -\varepsilon\phi(a_1)-\varepsilon\phi(\tilde a_1)\cr
   =&\phi(a)+\phi(\tilde a).\cr}$$
This shows that $\phi$ is a linear map from $A$ to $B$. Furthermore,
$T(a_1,a_2)=T_\phi(a_1,a_2)$ if $\{a_1,a_2\}$ spans a spacelike $2$ plane. We use Lemma
\brefe\ to see $T=T_\phi$. \qed\enddemo

\sethead\qctSE\head\S\qctSE\ The proof of Theorems \arefb\ (2) and \arefd\ (2)\endhead

\setref\erefa
\subhead\erefa\ The proof of Theorem \arefd\ (2) for metrics of
arbitrary signature\endsubhead
Let $q_A\ge5$. We use Lemma \brefb\ (2) to find $\psi$ which
relates
$\langle\cdot,\cdot\rangle_B$ to a positive definite metric $(\cdot,\cdot)_+$. Let $T$ be
admissible. Then $\psi T$ is admissible and takes values in
$\so(B,(\cdot,\cdot)_+)$. We use Lemma \drefc\ to see
\seteqn\erefaa
$$\eqalign{
&\psi T(a_1,a_2)b=(\hat\chi(a_1,a_2),b)_+\hat\xi
    -(\hat\xi,b)_+\hat\chi(a_1,a_2)\text{ or}\cr 
&\psi T(a_1,a_2)b=\pm\{(\hat\phi(a_2),b)_+\hat\phi(a_1)
   -(\hat\phi(a_1),b)_+\hat\phi(a_2)\}.}\tag\erefaa$$
Since $\psi^2=id$, we may apply $\psi$ to equation (\erefaa) and use
Lemma \brefb\ (1) to see:
$$\eqalign{
&T(a_1,a_2)b=\langle\psi\hat\chi(a_1,a_2),b\rangle_B\psi\hat\xi
    -\langle\psi\hat\xi,b\rangle_B\psi\hat\chi(a_1,a_2)\text{ or}\cr 
&T(a_1,a_2)b=\pm\{\langle\psi\hat\phi(a_2),b\rangle_B\psi\hat\phi(a_1)
   -\langle\psi\hat\phi(a_1),b\rangle_B\psi\hat\phi(a_2)\}.
\cr}$$
We complete the proof by setting $(\chi,\xi):=(\psi\hat\chi,\psi\hat\xi)$ or
$\phi:=\psi\hat\phi$ as appropriate. \qed

\setref\erefb
\subhead\erefb\ The proof of Theorem \arefb\ (2)\endsubhead Let
$\langle\cdot,\cdot\rangle$ be a non-degenerate metric of signature $(p,q)$ on $V$ where
$q\ge5$. Let
$R$ be an algebraic curvature tensor which has constant rank 2. We apply Theorem \arefd\ (2).
Suppose first that $R=R_{\chi,\xi}$.
Choose $\{x,y\}$ a spacelike orthonormal set so $\{\xi,x,y\}$ is an orthogonal set. Choose
$z$ so $\langle\xi,z\rangle\ne0$. We apply the Bianchi identity (\arefb) to see:
$$\eqalign{
  0=&\{\langle\chi(x,y),z\rangle+\langle\chi(y,z)x\rangle+\langle\chi(z,x),y\rangle\}\xi\cr
    &-\langle\xi,z\rangle\chi(x,y)-\langle\xi,x\rangle\chi(y,z)-\langle\xi,y\rangle
    \chi(z,x).}$$
Since $\langle\xi,z\rangle\ne0$ but $\langle\xi,x\rangle=\langle\xi,y\rangle=0$
we see $\chi(x,y)$ is a multiple of $\xi$ so 
$$\dim\{\range(R(x,y))\}\le1$$
which is false.
Thus $R=\pm R_\phi$ where $\ker(\phi)\cap V^+=\emptyset$. By replacing $R$ by $-R$ we may
suppose $R=R_\phi$. Since $R_\phi$ is an algebraic curvature tensor, we may apply Lemma
\brefc\ (2) to see $\phi=\phi^*$.
\qed

\sethead\qctSF\head\S\qctSF\ The classification of Ivanov-Petrova algebraic curvature
tensors of rank 2\endhead Let $R$ be a rank $2$ algebraic curvature tensor of rank $2$. We
use Theorem
\arefb\ to express $R=\pm R_\phi$ where $\phi$ is admissible. Let
$a_1\in\Bbb{S}^+(V)$. Choose $S\in\Cal{S}_3(V)$ so $a_1\in S$. Choose $a_2\in\Bbb{S}^+(S)$
so $a_2\perp a_1$ and so $\phi(a_2)\perp\phi(a_1)$. Let
$\pi:=\span\{a_1,a_2\}$ and let $\sigma:=\range(R_\phi(\pi))$. Since
$\{\phi(a_1),\phi(a_2)\}$ is an orthogonal basis for $\sigma$, we have
$$R_\phi(\pi):\phi(a_1)\rightarrow-\langle\phi(a_1),\phi(a_1)\rangle\phi(a_2)\text{ and }
  R_\phi(\pi):\phi(a_2)\rightarrow\langle\phi(a_2),\phi(a_2)\rangle\phi(a_1).$$
Let $\langle\cdot,\cdot\rangle_\sigma$ be the restriction of $\langle\cdot,\cdot\rangle$ to
$\sigma$. Since $\langle\phi(a_1),\phi(a_2)\rangle=0$, $\langle\cdot,\cdot\rangle_\sigma$ is
determined by the two inner products
$\langle\phi(a_1),\phi(a_1)\rangle$ and
$\langle\phi(a_2),\phi(a_2)\rangle$. We can relate the eigenvalue
structure of
$R_\phi(\pi)$ to $\langle\cdot,\cdot\rangle_\sigma$, there are 5 cases:
\roster
\smallskip\item The inner product $\langle\cdot,\cdot\rangle_\sigma$ is positive or
negative definite. Then
$R_\phi(\pi)$ has two non-zero purely imaginary eigenvalues $\pm\sqrt{-1}\lambda$ where
$\lambda^2=\langle\phi(a_1),\phi(a_1)\rangle\langle\phi(a_2),\phi(a_2)\rangle$.
\smallskip\item The inner product $\langle\cdot,\cdot\rangle_\sigma$ is indefinite
and non-degenerate.
Then
$R_\phi(\pi)$ has two non-zero real eigenvalues $\pm\lambda$ where
$\lambda^2=-\langle\phi(a_1),\phi(a_1)\rangle\langle\phi(a_2),\phi(a_2)\rangle$.
\smallskip\item The inner product $\langle\cdot,\cdot\rangle_\sigma$
is degenerate but non-trivial. Then
$R_\phi^2\ne0$ but $R_\phi^3=0$.
\smallskip\item The inner product $\langle\cdot,\cdot\rangle_\sigma$ is trivial. Then
$R_\phi^2=0$.
\endroster
\medbreak We can now prove one implication of Theorem \arefc. Suppose either that there
exists
$C\ne0$ so that $\langle\phi(v_1),\phi(v_2)\rangle=C\langle v_1,v_2\rangle$ for all
$v_i$ or that
$\langle\phi(v_1),\phi(v_2)\rangle=0$ for all
$v_i$. Then either case (1) or case (4) holds and $R$ is Ivanov-Petrova.

\medbreak To prove the other implication of Theorem \arefc, we suppose that $R_\phi$ is
Ivanov-Petrova. We first show that the complex Jordan normal form described by cases
(2) and (3) can not hold. Let $v_1\in\Bbb{S}(V^+)$. Since $q\ge5$, we may choose
$S\in\Cal{S}_5$ so that
$v_1\in S$. Choose $v_2\in S$ as above. Choose $v_3\in S$ so $v_3\perp v_i$ and
$\phi(v_3)\perp\phi(v_i)$ for $i=1,2$. Then $\{v_1,v_2,v_3\}$ is a spacelike orthonormal set
and $\{\phi(v_1),\phi(v_2),\phi(v_3)\}$ is an orthogonal linearly independent set. For $1\le
i<j\le 3$, let
$\pi_{ij}:=\span\{v_i,v_j\}$.
Suppose  $R_\phi(\pi_{ij})$ has two non-zero real eigenvalues. Then:
$$\langle\phi(v_i),\phi(v_i)\rangle\cdot\langle\phi(v_j),\phi(v_j)\rangle<0
\text{ for }1\le i<j\le 3.$$
This is not possible as at least two of the inner
products must have the same sign. Suppose $R_\phi(\pi_{ij})^3=0$ but
$R_\phi(\pi_{ij})^2\ne0$. Then the set of inner products
$$\{\langle\phi(v_i),\phi(v_i)\rangle,\langle\phi(v_j),\phi(v_j)\rangle\}$$
consists of one zero number and one non-zero number for $1\le i<j\le 3$. Again, this is not
possible as there are 3 such pairs.

We now deal with the remaining eigenvalue structures. Suppose $R_\phi(\pi)$ has two non-zero
purely imaginary eigenvalues $\pm\sqrt{-1}\lambda$ for all oriented spacelike $2$ planes
$\pi$. Let
$v_1$ and
$v_2$ be arbitrary unit spacelike vectors in
$V$. Since $q\ge5$, we can choose a third unit spacelike vector $v_3$ so $v_3\perp v_1$,
$v_3\perp v_2$, $\phi(v_3)\perp\phi(v_1)$, and so $\phi(v_3)\perp\phi(v_2)$. Let
$\pi_i:=\span\{v_i,v_3\}$ for $i=1,2$.
We then have
$$\eqalign{
 &0\ne\lambda^2=\langle\phi(v_1),\phi(v_1)\rangle
  \cdot\langle\phi(v_3),\phi(v_3)\rangle
 =\langle\phi(v_2),\phi(v_2)\rangle\cdot\langle\phi(v_3),\phi(v_3)\rangle\text{ so}\cr
  &\langle\phi(v_1),\phi(v_1)\rangle=\langle\phi(v_2),\phi(v_2)\rangle=C\text{ for all }
v_i\in\Bbb{S}(V^+)\text{ where }C=\pm\lambda.}$$
We rescale this identity to see that $\langle\phi(v),\phi(v)\rangle=C\langle v,v\rangle$
for all spacelike
$v$. Now let $v$
be arbitrary. Choose $v_1$ spacelike. Then
$v+\varepsilon v_1$ is spacelike for large $\varepsilon$ and hence 
$$\langle \phi(v+\varepsilon v_1),\phi(v+\varepsilon v_1)\rangle=C\langle v+\varepsilon
v_1,v+\varepsilon v_1\rangle\text{ for }\varepsilon>\varepsilon_0.$$
Since this equation is quadratic in the parameter $\varepsilon$,
it must hold $\varepsilon=0$ so
$$\langle\phi(v),\phi(v)\rangle=C\langle v,v\rangle\text{ for all }v\in V.$$

Suppose that $R_\phi(\pi)^2=0$ for any spacelike $2$ plane $\pi$. Then
$\range(R_\phi(\pi))$ is totally isotropic. Let
$v_1$ and
$v_2$ be arbitrary (not necessarily distinct) vectors in $V$. Choose unit spacelike vectors
$w_1$ and
$w_2$ so $w_1\perp w_2$. Let $v_i(\varepsilon):=v_i+\varepsilon w_i$. For large
$\varepsilon$, $\span\{v_1(\varepsilon),v_2(\varepsilon)\}$ is spacelike so
$\langle\phi(v_1(\varepsilon)),\phi(v_2(\varepsilon))\rangle=0$. Since this identity is
quadratic, it continues to hold for $\varepsilon=0$ and thus
$\langle\phi(v_1),\phi(v_2)\rangle=0$ for all $v_i\in V$ so $\range(\phi)$ is totally
isotropic. This completes the proof of Theorem
\arefc.

\sethead\qctSG\head\S\qctSG\ Realizing rank $2$ algebraic curvature tensors\endhead
Let $W$ be a vector space with a metric of signature $(r,s)$ and let $F:M\rightarrow W$ be an
immersion of a manifold $M$ into $W$ as a non-degenerate hypersurface; this means that the
first fundamental form is a non-degenerate metric of signature $(p,q)$ on $M$ where
$(p,q)=(r-1,s)$ or $(p,q)=(r,s-1)$. Choose a normal
$\nu$ along $M$ so $\langle\nu,\nu\rangle_W=\pm1$ and so $\langle
dF,\nu\rangle_W=0$. Let
$L$ be the second fundamental form and let $S$ be the associated shape operator determined
by the immersion. These two tensors are related by the identity:
$$L(v_1,v_2)=\langle v_1v_2(F),\nu\rangle_W=\langle Sv_1,v_2\rangle_W.$$
\setref\grefa
\proclaim{\grefa\ Lemma} If $(p,q)=(r-1,s)$, then $R=-R_S$; if $(p,q)=(r,s-1)$, then $R=R_S$.
\endproclaim

\demo{Proof} Although this is well known, we sketch the proof to establish later notation.
Fix a point $P\in M$. Choose a basis $\{e_1,...,e_{m+1}\}$ for $W$
so that
$\{e_1,...,e_m\}$ is a basis for $T_PM$ and so that $e_{m+1}=\nu(P)$. Let
$\{x^1,...,x^{m+1}\}$ be the dual basis for $W^*$; these give coordinates on $W$. Then $\vec
y:=(F^*x^1,...,F^*x^m)$ is a system of local coordinates on $M$ and the immersion takes the
form
$F(\vec y)=y^1e_1+...+y^me_m+f(\vec y)e_{m+1}$ where $f$ is a smooth scalar valued function
with
$df(P)=0$. Let $\varepsilon=\langle e_{m+1},e_{m+1}\rangle_W=\pm1$. We complete the
proof by computing:
$$\eqalign{
  &\langle S\partial_i,\partial_j\rangle(P)
   =L(\partial_i,\partial_j)(P)=\langle\partial_i\partial_jF,e_{m+1}\rangle_W(P)
     =\varepsilon_{m+1}\partial_i\partial_jf(P),\cr
  &g_{ij}=\langle e_i,e_j\rangle_W+\varepsilon\partial_i(f)\partial_j(f)
  ,\ \partial_k(g_{ij})(P)=0,\cr
  &R_{ijkl}(P)=\textstyle\frac12\{\partial_i\partial_kg_{jl}+\partial_j\partial_lg_{ik}
    -\partial_i\partial_lg_{jk}-\partial_j\partial_kg_{il}\}(P)\cr
  &\quad=\varepsilon\{\partial_i\partial_l(f)\cdot\partial_j\partial_k(f)
    -\partial_i\partial_k(f)\cdot\partial_j\partial_l(f)\}(P),\cr
  &(R_S)_{ijkl}(P)=\{\langle S\partial_i,\partial_l\rangle\langle S\partial_j,\partial_k
   \rangle-\langle S\partial_i,\partial_k\rangle\langle
S\partial_j\partial_l\rangle\}(P)=\varepsilon R_{ijkl}(P).\ \qed}$$
\enddemo

\demo{Proof of Theorem \arefe} By Theorem \arefb, we may choose $\phi$ so $R=\varepsilon
R_\phi$ where $\varepsilon=\pm1$; $\phi$ is symmetric and $\ker(\phi)$ contains no spacelike
vectors. We adopt the notation of equation (\arefdA) to define the germ of an immersion
$F_\phi$ with associated first fundamental form $g_\phi$, second fundamental form $L_\phi$,
and shape operator $S_\phi$. If $\{e_i\}$ is a basis for $V$ with associated coordinates
$\{y^i\}$, then
$$\eqalign{
   &F_\phi(\vec y):=y^1e_1+...+y^me_m+\textstyle\frac12
    \textstyle\sum_{i,j}y^iy^j\langle\phi(e_i),e_j\rangle_V\cdot
    e_{m+1},\cr
   &\partial_i(F_\phi)=e_i+
    \textstyle\sum_ky^k\langle\phi(e_i),e_k\rangle_V\cdot e_{m+1},\text{ and}\cr
   &g_{ij}(\vec y)=\langle e_i,e_j\rangle_V+\varepsilon\textstyle\sum_{k,l}
     y^ky^l\langle\phi(e_i),e_k\rangle_V\cdot
     \langle\phi(e_j),e_\ell\rangle_V.\cr}$$
This shows that $g_{ij}(0)=\langle\cdot,\cdot\rangle_V$. Since $e_{m+1}$ is the normal at
$\vec y=0$, we have that
$L_\phi(0)(e_i,e_j)=\varepsilon\langle\phi(e_i),e_j\rangle_V$. Thus
$S_\phi(0)=\varepsilon\phi$ so by Lemma \grefa,
${}^gR=\varepsilon R_{\varepsilon\phi}=R$. Consequently $(V,g_\phi)$ is a
geometric realization of
$R$.

We may choose the basis $\{e_i\}$ so that $\{e_i\}_{i>s}$ is a basis for $\ker(\phi)$;
if $\phi$ is invertible, then $s=m$ and this is the empty basis. If $i>s$, we have
$\partial_i\partial_j(F)=\langle\phi(e_i),e_j\rangle_V=0$ so
$L_\phi(P)(\partial_i,\partial_j)=0$ and $S_\phi(P)(\partial_i)=0$ at any point $P$. Since
we are only considering an arbitrarily small neighborhood of the origin,
$\rank(S_\phi(P))\ge\rank(S_\phi(0))=s$. Thus equality holds and
$\ker(S_\phi(P))=\span\{e_i\}_{i>s}$. If $i,j>s$, then $g_{ij}(P)=g_{ij}(0)$. Since
$\ker(\phi)=\span\{e_i\}_{i>s}$ contains no spacelike vectors, $\ker(S_\phi)(P)$ contains no
spacelike vectors and $(V,g_\phi)$ is the germ of a rank $2$ pseudo-Riemannian manifold. 
\qed\enddemo

\sethead\qctSH\head\S\qctSH\ Examples\endhead
\setref\hrefa
\subhead\hrefa\ Example\endsubhead We have assumed in the previous sections that $q_A\ge5$.
We now show that Theorem
\arefd\ (2) fails if $q_A=4$. Let $\{e_1^-,...,e_p^-,e_1^+,...,e_q^+\}$ be a basis
for $\Bbb{R}^{(p,q)}$. We define $\langle\cdot,\cdot\rangle_{(p,q)}$ by:
$$\langle e_i^-,e_j^-\rangle_{(p,q)}=-\delta_{ij},\ \langle e_i^-,e_j^+\rangle_{(p,q)}=0,
  \text{ and }\langle e_i^+,e_j^+\rangle_{(p,q)}=\delta_{ij}.$$
We suppose first $p_A=0$. Let $\text{orn}:=e_1^+\wedge e_2^+\wedge e_3^+\wedge e_4^+$ be the
standard orientation of
$\Bbb{R}^{(0,4)}$. Extend the inner product on $\Bbb{R}^{(0,4)}$ to an
inner product on
$\Lambda^4(\Bbb{R}^{(0,4)})$ and define: 
$$\eqalign{
 &T^{(0,4)}:\Bbb{R}^{(0,4)}\otimes\Bbb{R}^{(0,4)}\rightarrow
     \so(\Bbb{R}^{(0,4)},\langle\cdot,\cdot\rangle_{(0,4)})\text{ by}\cr
 &\langle T^{(0,4)}(v_1,v_2)v_3,v_4\rangle_{(0,4)}=\langle v_1\wedge v_2\wedge v_3\wedge v_4,
\text{orn}\rangle_{(0,4)}.}$$
Let $\{f_1,f_2\}$ be an oriented orthonormal basis for a $2$ plane $\pi$. We
extend $\{f_1,f_2\}$ to an oriented orthonormal basis $\{f_1,f_2,f_3,f_4\}$ for
$\Bbb{R}^{(0,4)}$ and compute:
$$\eqalign{
  &T^{(0,4)}(f_1,f_2)f_1=0,\phantom{ and }T^{(0,4)}(f_1,f_2)f_2=0,\cr
  &T^{(0,4)}(f_1,f_2)f_3=f_4,\text{ and }
  T^{(0,4)}(f_1,f_2)f_4=-f_3.}$$
Thus $T^{(0,4)}$ has constant rank $2$. Since $\range(T^{(0,4)}(\pi))=\pi^\perp$,
$\Phi([a])=\{0\}$. Thus $T^{(0,4)}\ne T^{(0,4)}_\phi$ and $T^{(0,4)}\ne
T^{(0,4)}_{\chi,\xi}$ for any $\phi$ or
$(\chi,\xi)$. More generally, let $\sigma:\Bbb{R}^{(p,4)}\rightarrow\Bbb{R}^{(0,4)}$ denote
projection on the last $4$ components. We extend $T^{(0,4)}$ to
$\Bbb{R}^{(p,4)}\otimes\Bbb{R}^{(p,4)}$ and thereby construct an algebraic
curvature tensor which does not satisfy the conclusions of Theorem \arefd\ (2)
for any
$p$ by defining:
$$T^{(p,4)}(v_1,v_2):=T^{(0,4)}(\sigma(v_1),\sigma(v_2)).$$

\setref\hrefb
\subhead\hrefb\ Example\endsubhead Let $q=p+1\ge5$. We define
$$\phi(e_i^\pm)=\pm(e_i^-+e_i^+)\text{ for }i\le p\text{ and }
  \phi(e_{p+1}^+)=e_{p+1}^+.$$
Then $\phi$ is self-adjoint and $\ker(\phi)=\span\{e_i^--e_i^+\}_{i\le p}$ contains no
spacelike vectors so $R_\phi$ has constant rank $2$. Let $\pi$ be a spacelike $2$ plane.
Since the metric on
$\range(R(\pi))$ is degenerate, $R_\phi^3(\pi)=0$ so $0$ is the only eigenvalue of $R_\phi$.
Let
$\pi_i:=\span\{e_1^+,e_i^+\}$. Then
$R_\phi(\pi_2)^2=0$ and $R_\phi(\pi_{p+1})^2\ne0$ so
$R_\phi$ is not an Ivanov-Petrova algebraic curvature tensor. When dealing with an
indefinite metric, it is necessary to study the complex Jordan normal form and not just
the eigenvalue structure.

\setref\hrefc
\subhead\hrefc\ Example\endsubhead We now construct an algebraic curvature tensor which is
(spacelike) rank $2$ Ivanov-Petrova but which does not have constant rank on the Grassmanian
of oriented timelike
$2$ planes as follows. Let
$q=p-1\ge5$. We define:
$$\phi(e_i^\pm)=\pm(e_i^-+e_i^+)\text{ for }i\le q,\text{ and }
  \phi(e_{q+1}^-)=0.$$
Since $\phi$ is self-adjoint and 
$\ker(\phi)=\span\{e_i^-+e_i^+,e_{q+1}^-\}_{i\le q}$ contains no
spacelike vectors, $R_\phi(\pi)$ has rank $2$. Furthermore as
the range of $\phi$ is totally isotropic, $R_\phi$ is (spacelike) rank $2$ Ivanov-Petrova.
Let $\pi_i^-:=\span\{e_1^-,e_i^-\}$ be oriented timelike $2$ planes. Then
$\rank(R_\phi(\pi_i))=2$ if $2\le i\le q$ while $\rank(R_\phi(\pi_i))=0$ if $i=q+1$. Thus
the rank of $R_\phi$ is not constant on the timelike $2$ planes.

\subhead Acknowledgements\endsubhead 
The research of Gilkey was partially supported by the NSF (USA) and MPI (Leipzig, Germany).
It is a pleasant task to thank Professors Cortes, Simon, and Vrancken for helpful
discussions concerning the material of \S\qctSG.

\medbreak\noindent PG: Mathematics Department, University of Oregon, Eugene Or 97403 USA
\par\noindent\phantom{PG: }{\it E-mail address}: gilkey\@darkwing.uoregon.edu
\medbreak\noindent TZ: Mathematics Department, Murray State University, Murray Kentucky
42071 USA
\par\noindent\phantom{TZ: }{\it E-mail address}: tan.zhang\@murraystate.edu
\footnote" "{Version \version\ printed\ \number \day\ \nmonth\ \number\year}
\enddocument
\bye